\documentstyle[12pt]{article}
\catcode `\@=11
\@addtoreset{equation}{section}

\newcommand{\be}{\begin{equation}}
\newcommand{\ee}{\end{equation}}
\newcommand{\R}{\rm I \mkern -3mu R}
\newcommand{\N}{\rm I \mkern -3mu N}
\newcommand{\Hil}{\cal H}
\newcommand{\C}{\cal C}
\newcommand{\K}{\cal K}
\newcommand{\U}{\cal U}
\newcommand{\W}{\cal W}
\newcommand{\Gpo}{{\cal G}^{\uparrow}_{+}}
\newcommand{\Go}{{\cal G}^{\uparrow}}
\newcommand{\Gpocov}{\widetilde{{\cal G}^{\uparrow}_{+}}}
\newtheorem{thm}{Theorem}
\newtheorem{rem}{Remark}
\newtheorem{lemma}{Lemma}
\newtheorem{prop}{Proposition}
\begin{document}
\thispagestyle{empty}
\begin{flushright}
IFA-FT-403-1994, December
\end{flushright}
\bigskip\bigskip\begin{center}
{\bf \Large{PROJECTIVE SYSTEMS OF IMPRIMITIVITY}}
\end{center}
\vskip 1.0truecm
\centerline{\bf
D. R. Grigore\footnote{e-mail: grigore@roifa.bitnet, grigore@ifa.ro}}
\vskip5mm
\centerline{Dept. Theor. Phys., Inst. Atomic Phys.,}
\centerline{Bucharest-M\u agurele, P. O. Box MG 6, ROM\^ANIA}
\vskip 2cm
\bigskip \nopagebreak \begin{abstract}
\noindent
We apply Mackey procedure of classifying projective systems of
imprimitivity to a thorough study of the projective unitary
irreducible representations of the Galilei group in 1+3 and 1+2 dimensions.
\end{abstract}
\newpage\setcounter{page}1

\section{Introduction}

Recently we have studied the projective unitary irreducible
representations of the Poincar\'e group in 1+2 dimensions
\cite{G1}. Here we provide the same analysis for the Galilei
group in 1+2 dimensions.

The main difficulty in applying straighforwardly Mackey method
is due to the rather complicated structure of the multiplier
group (see \cite{MS} Appendix A, \cite{G2}). One can proceed in
two ways. The first one mimicks the usual method of Mackey as
presented, for instance in \cite{V}, and is based on the use of
group extensions. This line of argument has been developped in
\cite{G2}. The second method is based on the use of projective
systems of imprimitivity and is also due to Mackey \cite{M}. We
will pursue the second method in this paper, but for the benefit
of the reader, we will also outline the study of the projective
systems of imprimitivity. We will prefer to follow the general
line of argument and the notations of \cite{V} and indicate the
necessary changes. This will be done in Section 2. In Section 3
we apply the method to the Galilei group in 1+3 (only for
illustrative purpose) and afterwards to the case of 1+2
dimensions.

\section{Projective systems of imprimitivity}

2.1 Let
$G$
be a Borel group acting in the Borel space
$X$
and let
$\Hil$
be a Hilbert space. A {\it projective system of imprimitivity for
$G$
based on
$X$
and acting in}
$\Hil$
is a pair
$
(U,P)
$
where
$
P~(E \mapsto P_{E})
$
is a projector valued measure based on
$X$
and acting in
$\Hil$,
$U~(g \mapsto U_{g})$
is a projective representation of
$G$
in
$\Hil$
and such that the relations
\be
U_{g} P_{E} U_{g}^{-1} = P_{g\cdot E}
\ee
are satisfied for all
$
g \in G
$
and all Borel sets
$
E \subset X.
$
If
$U$
is a
$m$-representation i.e.
$
m: G \times G \rightarrow {\bf T}
$
(here {\bf T} is the set of complex numbers of modulus 1) is the
multiplier of
$U$:
\be
U_{g_{1}g_{2}} = m(g_{1},g_{2})~U_{g_{1}}~U_{g_{2}}
\ee
for all
$
g_{1}, g_{2} \in G,
$
then one also says that
$
(U,P)
$
is a {\it m-system of imprimitivity}. Now the notions of {\it
equivalence}, {\it irreducibility}, {\it direct sum} and {\it
commuting ring} for a
$m$-system
of imprimitivity are defined as for the usual case
$
m = 1
$
(see \cite{V}, ch. VI, pp. 203).

2.2 If one wants to classify
$m$-systems
of imprimitivity the first step is to connect them with
$m$-cocycles. Let
$G$
be a Borel group acting in the Borel space
$X$
and let
$\C$
be an invariant measure class on
$X$.
Suppose
$M$
is a standard Borel group with the identity denoted by {\bf 1}.
(The usual case will be
$M = {\cal U}(\Hil) \equiv$
the group of all the unitary operators in the Hilbert space
$\Hil$).
Finally, let
$
m: G \times G \rightarrow {\bf T}
$
be a multiplier of
$G$.
We say that
$
f: G \times X \rightarrow M
$
is a
$m-(G,X,M)$-{\it cocycle}
(or shortly a
$m$-cocycle) if

(i)
$f$
is a Borel map

(ii)
$
f(e,x) = 1
$
for almost all
$
x \in X
$

(iii)
$
f(g_{1}g_{2},x) = m(g_{1},g_{2})~~f(g_{1},g_{2}\cdot x)~~f(g_{2},x)
$
for almost all
$
(g_{1},g_{2},x) \in G \times G \times X.
$

The notions of {\it strict cocycle}, {\it (strict) cohomology},
{\it (strict) coboundary}, {\it (strict) cohomology class} are
defined as for the usual case
$
m = 1.
$
We denote, as usual, (strict) equivalence by
$
(\approx) \sim
$.

2.3 The connection between
$m$-systems
of imprimitivity and
$m$-cocycles
goes practically unchanged, as in \cite{V}. So, lemma
6.6 from \cite{V} trivially modifies in:

\begin{lemma}
Let
$G$
be a Borel group verifying the second axiom of countability and
$\Hil$
a separable Hilbert space. Suppose that
$
L~(g \rightarrow L_{g})
$
is a mapping of
$G$
into the set of bounded operators in
$\Hil$
such that
(i) for all
$
f,f' \in {\cal H},~ g \mapsto <f,L_{g}f'>
$
is a Borel map; (ii)
$
L_{g}
$
is unitary for almost all
$g$;
(iii)
$
L_{g_{1}g_{2}} = m(g_{1},g_{2}) L_{g_{1}}~L_{g_{2}}
$
almost everywhere in
$
G \times G.
$

Then there exists exactly one m-representation
$U$
of
$G$
in
$\Hil$
such that
$
L_{g} = U_{g}
$
for almost all
$g$.
\end{lemma}

Theorem 6.7 of \cite{V} goes into:

\begin{thm}
Let the notations be as above. Let
$\K$
be a Hilbert space,
$
M = {\cal U}(\K),
$
$X$
a
$G$-Borel
space,
$\C$
a
$G$-invariant
measure class on
$X$
and
$
\alpha \in \C
$
a measure. For each
$
g \in G,
$
let
$
r_{g}
$
be a Borel function which is a version of the Radon-Nycodim
derivative
$
d\alpha/d\alpha^{g^{-1}}.
$
Suppose that
$\phi$
is a
$m-(G,X,M)$-cocycle relative to the measure class of
$\alpha$.
Then, there exists a unique m-system of imprimitivity
$
(U,P)
$
acting in
$
{\cal H} \equiv L^{2}(X,{\cal K},\alpha)
$
such that:
\be
P_{E} f = \chi_{E}~f,~~\forall f \in \Hil
\label{proj}
\ee
and for almost all
$
g \in G
$
we have
\be
\left(U_{g}~f\right)(x) = \left\{ r_{g}(g^{-1}\cdot
x)\right\}^{1/2}~\phi(g,g^{-1}\cdot x)~f(g^{-1}\cdot x)
\label{repr}
\ee
for all
$
f \in \Hil
$
and for almost all
$
x \in X.$
Moreover, the equivalence class of
$
(U,P)
$
depends only on the measure class
$\C$
and on the cohomology class of
$\phi.$
\label{co-sys}
\end{thm}

2.4 Conversely, going from
$m$-sytems
of imprimitivity to
$m$-cocycles
involves again only a slight departure from \cite{V}.
So, lemma 6.10 of \cite{V} stays true:

\begin{lemma}
If the system of imprimitivity
$
(U,P)
$
is irreducible, then the measure class
$\C$
is ergodic. If
$\C$
is ergodic, then
$P$
is homogeneous.
\label{ergo}
\end{lemma}

Next, suppose that
$
P = P(\{ {\cal K}_{n} \},\{\alpha_{n} \})
$
is the decomposition of
$P$
into homogeneous projection valued measures. So
$
{\cal K}_{n}
$
is a
$n$-dimensional
Hilbert space,
$
n \in \N \cup \{\infty\},~\alpha_{1}, \alpha_{2},...,\alpha_{\infty}
$
are mutually singular
$\sigma$-finite
measures on
$X$
and
$
\alpha = \sum_{n=1}^{\infty} \alpha_{n}.
$
In the Hilbert space
$
{\cal H}_{n,\alpha} = L^{2}(X,{\cal K}_{n},\alpha)
$
we can construct a
$m$-system
of imprimitivity if we have at our disposal a
$m-(G,X,M)$-cocycle;
indeed we simply use Theorem \ref{co-sys} and define
$P$
and
$U$
according to (\ref{proj}) and (\ref{repr}) respectively. We denote
in this case
$P$
by
$
P^{n,\alpha}.
$

Then we have the counterpart of Theorem 6.11 of \cite{V}:

\begin{thm}
Let
$
(U,P)
$
be a m-system of imprimitivity acting in the Hilbert space
$\Hil$
such that
$P$
is homogoneous of multiplicity
$
n~(1 \leq n \leq \infty).
$
Let
$\alpha$
be a
$\sigma$-finite
measure in the measure class
$\C$
of
$P$.
Then
$
(U,P)
$
is equivalent to a m-system of imprimitivity
$
(U',P^{n,\alpha})
$
acting in
$
{\cal H}_{n,\alpha}.
$
Moreover, there exists a one-one correspondence between the set
of all cohomology classes of
$
m-(G,X,{\cal U}({\cal K}_{n}))$-cocycles relative to
$
\C_{P}
$
and the set of all equivalence classes of m-systems of
imprimitivity of the form
$
(U',P^{n,\alpha}).
$
\end{thm}

2.5 From the preceeding two subsections it is clear that one is
reduced to the classification problem for
$
m-(G,X,M))$-cocycles.
As in \cite{M}, \cite{V} this can be done rather completely in
the case of a transitive action of
$G$
on
$X$.
In this case, one chooses
$
x_{0} \in X
$
arbitrary; it is known that the stability subgroup
$
G_{0} \equiv G_{x_{0}}
$
is closed and we have
$
X \simeq G/G_{0}
$
as Borel spaces. So, from now on
$
X = G/G_{0}
$
for some closed subgroup
$
G_{0}.
$
We denote by
$
\beta: G \rightarrow X
$
the canonical projection:
$
\beta(g) = g\cdot G_{0}$;
it is known that
$\beta$
is a Borel map. Next one notices that if
$
\alpha_{0}
$
is a
$\sigma$-finite
$G$-quasi-invariant
measure on
$G$,
then one can produce a
$\sigma$-finite
$G$-quasi-invariant
measure on
$X$
according to the formula:
\be
\alpha(A) = \alpha_{0}(\beta^{-1}(A))
\label{measure}
\ee
for all Borel subsets
$
A \subset X.
$

It is known that on a homogeneous
$G$-space
(as
$
X = G/G_{0}
$)
there exists a unique
$G$-invariant
measure class. So, if we take
$
\alpha_{0}
$
to be the Haar measure,
$\alpha$
from (\ref{measure}) will give us a representative from this
measure class. So, from now on, when speaking of cocycles on
$X$
we will always mean that they are relative to this measure class.

2.6 We concentrate for the moment on strict
$m-(G,X,M)$-cocycles.
Let
$X$
be a transitive
$G$-space
and
$f$
a strict
$m-(G,X,M)$-cocycle.
If
$x \in X$,
let us define
$
D: G_{x} \rightarrow M
$
by:
\be
D(h) \equiv f(h,x).
\ee

Then
$D$
is a
$m$-representation of the group
$
G_{x}
$
in
$M$.
One calls
$D$
the {\it m-homomorphism defined by f in}
$
x \in X$.
Then we have two elementary results which generalize lemma 5.23
of \cite{V}.

\begin{lemma}
Let
$f$
be a (strict)
$m-(G,X,M)$-cocycle.
We define
$
f^{0}: G \times G \rightarrow M
$
by:
\be
f^{0}(g,g') \equiv f(g,\beta(g')),~~\forall g,g' \in G.
\ee
Then
$
f^{0}
$
is a (strict)
$m-(G,G,M)$-cocycle.
Conversely, let
$F$
be a (strict)
$m-(G,G,M)$-cocycle,
$
G_{0} \subset G
$
a closed subgroup and
$
X = G/G_{0}.
$
Then there exists a (strict)
$m-(G,X,M)$-cocycle
$f$
such that
$
f^{0} = F
$
{\it iff}:
\be
F(g,g'h) = F(g,g'),~~\forall g,g' \in G,~\forall h \in G_{0}.
\label{cond}
\ee
Moreover this correspondence
$
f \leftrightarrow F
$
is one-one.
\label{f-f0}
\end{lemma}

\begin{lemma}
If
$F$
ia a strict
$
m-(G,G,M)$-cocycle,
there exists a unique Borel map
$
b: G \rightarrow M
$
such that:
\be
F(g,g') = m(g,g')^{-1}~b(g,g')~b(g')^{-1}.
\ee
\end{lemma}

Next, we come to a succesion of results generalizing lemma 5.24
of \cite{V}.

\begin{lemma}
A strict
$
m-(G,G,M)$-cocycle
$F$
verifies (\ref{cond}) for some closed subgroup
$
G_{0} \subset G
$
{\it iff} the map
$
b: G \rightarrow M
$
associated as above verifies:
\be
b(e) = {\bf 1}
\label{id}
\ee
\be
b(gh) = m(g,h)~b(g)~b(h).
\label{b}
\ee

In this case the map
$
D: G \rightarrow M
$
defined by
\be
D(h) \equiv b(h)
\ee
is a
$m$-representation of
$G_{0}$
in
$M$.
\end{lemma}

\begin{lemma}
In the conditions above let
$f$
be a strict
$
m-(G,X,M)$-cocycle,
$
f^{0}
$
the associated strict
$
m-(G,G,M)$-cocycle
and
$
D: G_{0} \rightarrow M
$
the map corresponding to
$
f^{0}
$
as above. Then
$f$
defines
$D$
in
$
G_{0} \in G/G_{0}.
$
\end{lemma}

\begin{lemma}
There is a one-one correspondence between strict
$
m-(G,X,M)$-cocycles
and maps
$
b: G \rightarrow M
$
verifying (\ref{id}) and (\ref{b}) above.
\label{one-one}
\end{lemma}

One can construct maps verifying the conditions above using
Borel cross sections i. e. Borel maps
$
c: X \rightarrow G
$
verifying
\be
\beta \circ c = id.
\ee

\begin{lemma}
Let
$
G_{0} \subset G
$
be a closed subgroup and
$
D:  G_{0} \rightarrow M
$
a m-representation. Then there exists a Borel map
$
b: G \rightarrow M
$
such that
$
b\vert_{G_{0}} = D
$
and
\be
b(gh) = m(g,h)~b(g)~b(h),~~\forall g \in G,~\forall h \in G_{0}.
\ee
\label{Db}
\end{lemma}

{\bf Proof}

Indeed, one knows that cross sections
$
c: X \rightarrow G
$
do exist. One first arranges such that:
\be
c(G_{0}) = e.
\ee

Then, one defines
$
a:G \rightarrow G
$
according to:
\be
a(g) \equiv c(\beta(g))^{-1} g.
\ee

In fact
$
a: G \rightarrow G_{0}.
$
If
$
D: G_{0} \rightarrow M
$
is a
$m$-representation one finally checks that the map
\be
b(g) \equiv m(c(\beta(g))^{-1},g)^{-1}~D(a(g))
\ee
verifies the conditions in the statement of the lemma.
$\Box$

\begin{lemma}
Let
$
f_{1}, f_{2}
$
be two strict
$
m-(G,X,M)$-cocycles
and
$
D_{1}, D_{2}
$
the associated m-representations. Then
$f_{1}$
and
$f_{2}$
are cohomologous
$
(f_{1} \approx f_{2})
$
{\it iff}
$D_{1}$
and
$D_{2}$
are equivalent m-representations.
\end{lemma}

\begin{lemma}
Let
$
D: G_{0} \rightarrow M
$
be a m-representation and
$
b_{1}, b_{2}:  G \rightarrow M
$
two Borel maps verifying
$
b_{i}\vert_{G_{0}} = D
$
and
$
b_{i}(gh) = b_{i}(g)~b_{i}(h),~~\forall g \in G,~\forall h \in
G_{0}~(i = 1,2).
$
Let
$
f_{1}, f_{2}
$
be the strict
$
m-(G,X,M)$-cocycles associated to
$
b_{1}, b_{2}
$
respectively (see lemma \ref{one-one}). Then
$
f_{1} \approx f_{2}.
$
\end{lemma}

In conclusion we have:
\begin{prop}
There exists a one-one correspondence between the set of
cohomology classes of strict
$
m-(G,X,M)$-cocycles
and the set of equivalence classes of
$
m\vert_{G_{0}\times G_{0}}$-representations of
$G_{0}$ in
$M$.
More precisely, let
$
D: G_{0} \rightarrow M
$
be a
$
m\vert_{G_{0}\times G_{0}}$-representation.
Then every strict
$
m-(G,X,M)$-cocycle
which defines
$D$
in
$G_{0}$
is strictly cohomologous to:
\be
\phi(g,x) = m(c(g\cdot x)^{-1},g)^{-1}~
m(c(g\cdot x)^{-1}gc(x),c(x)^{-1}) D(c(g\cdot x)^{-1}gc(x)).
\label{co}
\ee
\end{prop}

Needless to say, all the computations involved in proving the
results above are elementary.

2.7 We analyse here the general case of
$
m-(G,X,M)$-cocycles.
The key point is to generalize lemma 5.26 of \cite{V}. We do
this in some detail.

\begin{lemma}
If
$f$
is any
$
m-(G,X,M)$-cocycle
there exists a strict
$
m-(G,X,M)$-cocycle
$
f_{1}
$
such that
\be
f(g,x) = f_{1}(g,x)
\ee
for almost all
$
(g,x) \in G \times X;~f_{1}
$
is uniquely determined up to strict cohomology.
\end{lemma}

{\bf Proof} Follows \cite{V}.

If
$f$
is a
$
m-(G,X,M)$-cocycle,
let
$
f^{0}
$
be the associated
$
m-(G,G,M)$-cocycle
(see lemma \ref{f-f0}). We have
$$
f^{0}(g_{1}g_{2},g_{3}) = m(g_{1},g_{2})~f^{0}(g_{1},g_{2}g_{3})~
f^{0}(g_{2},g_{3})
$$
for almost all
$
(g_{1},g_{2},g_{3}) \in G \times G \times G.
$
Hence there exists a
$
g_{0} \in G_{0}
$
such that
$$
f^{0}(g_{1}g_{2},g_{0}) = m(g_{1},g_{2})~f^{0}(g_{1},g_{2}g_{0})~
f^{0}(g_{2},g_{0})
$$
for almost all
$
(g_{1},g_{2}) \in G \times G.
$
We take
$
g = g_{1},~g' = g_{2}g_{0}
$
and conclude that
$$
f^{0}(g,g') = m(g,g'g_{0}^{-1})^{-1}~f^{0}(gg'g_{0}^{-1},g_{0})~
f^{0}(g'g_{0}^{-1},g_{0})^{-1}
$$
or, using the multiplier identity:
$$
f^{0}(g,g') = m(g,g')^{-1}~m(g',g_{0}^{-1})~m(gg',g_{0}^{-1})~
f^{0}(gg'g_{0}^{-1},g_{0})~f^{0}(g'g_{0}^{-1},g_{0})^{-1}
$$
for almost all
$
(g,g') \in G \times G.
$

We define the map
$
d_{0}: G \rightarrow M
$
by
$$
d_{0}(g) \equiv m(g,g')^{-1}~f^{0}(gg_{0}^{-1},g_{0}).
$$

Then
$
d_{0}
$
is a Borel map and we have:
$$
f^{0}(g,g') = m(g,g')^{-1}~d_{0}(gg')~d_{0}(g')^{-1}
$$
for almost all
$
(g,g') \in G \times G.
$
If we redefine
$
d_{0}:
$
$$
d(g) \equiv d_{0}(g)~d_{0}(e)^{-1}
$$
then we have:
$$
d(e) = {\bf 1}
$$
and
\be
f^{0}(g,g') = m(g,g')^{-1}~d(gg')~d(g')^{-1}
\label{fd}
\ee
for almost all
$
(g,g') \in G \times G.
$

Now let
$
h \in G_{0}
$
be arbitrary. We have then:
$$
f^{0}(g,g'h) = f^{0}(g,g').
$$

Inserting here (\ref{fd}) we obtain:
$$
m(gg',h)^{-1}~d(gg')^{-1}~d(gg'h) = m(g',h)^{-1}~d(g')^{-1}~d(g'h)
$$
for almost all
$
(g,g') \in G \times G.
$
Because the measure class
$\C$
on
$X = G$
is known to be ergodic, it follows that the function:
$
g' \mapsto m(g',h)^{-1}~d(g')^{-1}~d(g'h)
$
is almost everywere constant i.e. we have
$
D: G_{0} \rightarrow M
$
such that;
\be
d(gh) = m(g,h)~d(g)~D(h)
\label{dD}
\ee
for almost all
$
g \in G
$;
this fixes
$D$
uniquely. One easily finds out that;
$$
D(e) = {\bf 1}
$$
$$
D(h_{1}h_{2}) = m(h_{1},h_{2})~D(h_{1})~D(h_{2})
$$
i.e
$D$
is a
$
m\vert_{G_{0} \times G_{0}}$-representation of
$G_{0}$
in
$M$.

We prove now that
$D$
is a Borel map. Indeed, let
$\lambda$
be a quasi-invariant measure on
$G$
normalized by
$
\lambda(G) = 1.
$
Then from (\ref{dD}) one has:
$$
D(h) = \int d\lambda(g)~m(g,h)^{-1}~d(g)^{-1}~d(gh).
$$

Applying lemma \ref{Db} one can find out a Borel map
$
d': G \rightarrow M
$
such that:
$$
d'(gh) = m(g,h)~d'(g)~D(h).
$$

Comparing with (\ref{dD}) it follows that
$
\forall h \in G_{0}:
$
$$
d(g)~d'(g)^{-1} = d(gh)~d'(gh)
$$
for almost all
$
g \in G.
$
Like in \cite{V} it follows that there exists a Borel map
$
k: X \rightarrow M
$
such that
$$
d(g) = k(\beta(g))~d'(g)
$$
for almost all
$
g \in G.
$
If we define
$
d_{1}: G \rightarrow M
$
by:
$$
d_{1} \equiv k(\beta(g))~d'(g)
$$
then we obviously have:
$$
d(g) = d_{1}(g)
$$
for almost all
$
g \in G
$
and also:
$$
d_{1}(gh) = m(g,h)~d_{1}(g)~D(h),~~~\forall g \in G,~\forall h
\in G_{0}
$$
so, according to lemma \ref{one-one} there exists a strict
$
m-(G,X,M)$-cocycle
$f_{1}$
such that
$$
f^{0}_{1}(g,g') = m(g,g')^{-1}~d_{1}(gg')~d_{1}(g')^{-1}.
$$

It is clear that
$
f = f_{1}
$
almost everywere in
$G \times X.$
$\Box$

As corollaries we have as in \cite{V}

\begin{lemma}
Let
$
f_{i}~(i = 1,2)
$
be two strict
$
m-(G,X,M)$-cocycles.
Then
$
f_{1} \sim f_{2} \Leftrightarrow f_{1} \approx f_{2}.
$
\end{lemma}

\begin{lemma}
Let
$\gamma$
be a cohomology class of
$
m-(G,X,M)$-cocycles
and
$
\gamma_{st} \subset \gamma
$
the set of strict
$
m-(G,X,M)$-cocycles.
Then
$
\gamma_{st}
$
is non-void and is a strict class of cohomology.
\end{lemma}

\begin{lemma}
There is a one-one correspondence between the set of equivalence
classes of transitive m-systems of imprimitivity for
$G$
in
$
{\cal H}_{n,\alpha}
$
and the set of strict cohomology classes of
$
m-(G,X,{\cal U}({\cal K}_{n}))$-cocycles
relative to
$
\C_{P}.
$
\end{lemma}

2.8 In conclusion, we formulate the main results.

\begin{prop}
There is a one-one correspondence between the set of equivalence
classes of transitive m-systems of imprimitivity for
$G$
based on
$
X
$
and the set of equivalence classes of
$
m\vert_{G_{0} \times G_{0}}$-representations
of
$G_{0}$
in
$M.$
\end{prop}

\begin{prop}
Let
$
D: G_{0} \rightarrow M
$
be a m-representation and
$
(U,P)
$
a m-system of imprimitivity associated to
$D$
in the sense of the preceeding proposition. Then the commuting
ring of
$D$
is isomorphic to the commuting ring of
$
(U,P).
$
Hence, the correspondence from the preceeding proposition
preserves irreducibility and the direct sum.
\end{prop}

The proofs of these propositions are similar to those in \cite{V}.
One has to notice that all the multiplier factors in various formul\ae~
conveniently cancel.

So, the procedure of constructing irreducible m-systems of
imprimitivity consists in following three steps.

(A) First, one notes that from lemma \ref{ergo} it easily
follows that the
$m$-system
of imprimitivity is transitive. So, one classifies all
$G$-orbits
in
$X$.
Let us fix now a certain orbit
${\cal O}$.
We have to classify the irreducible
$m$-systems
of imprimitivity for
$G$
based on
$
X = {\cal O}.
$
One identifies some closed subgroup
$G_{0}$
such that
$
X \simeq G/G_{0}
$
and lists the irreducible
$
m\vert_{G_{0} \times G_{0}}$-representations
up to unitary equivalence.

(B) One computes an associated
$
m-(G,X,M)$-cocycle using, eventually, the formula (\ref{co}).

(C) One constructs the
$m$-system
of imprimitivity according to (\ref{proj}) and (\ref{repr})
(see theorem \ref{co-sys}).

\section{The Projective Unitary Irreducible Representations
of the Galilei Group}

3.1 By definition the ortochronous Galilei group in
$1+n$
dimensions
$\Go$
is set-theoretically
$
O(n-1) \times \R^{n-1} \times \R \times \R^{n-1}
$
with the composition law:
\be
(R_{1},{\bf v}_{1},\eta_{1},{\bf a}_{1}) \cdot
(R_{2},{\bf v}_{2},\eta_{2},{\bf a}_{2}) =
(R_{1}R_{2},{\bf v}_{1}+R_{1}{\bf v}_{2},\eta_{1}+\eta_{2},
{\bf a}_{1}+R_{1}{\bf a}_{2}+\eta_{2}{\bf v}_{1}).
\ee

We organize
$
\R^{n-1}
$
as column vectors,
$
O(n-1)
$
as
$
(n-1) \times (n-1)
$
real orthogonal matrices and we use consistently matrix
notations. This group acts naturally on
$
\R \times \R^{n-1}
$
as follows
\be
(R,{\bf v},\eta,{\bf a}) \cdot(T,{\bf X}) =
(T+\eta,R{\bf X}+T{\bf v}+{\bf a}).
\ee

The proper orthochronous Galilei group
$\Gpo$
is by definition:
\be
\Gpo \equiv \{(R,{\bf v},\eta,{\bf a}) \vert det(R) = 1 \}.
\ee

The groups
$\Go$
and
$\Gpo$
are Lie groups.

As it is well known, the classification of all projective
unitary irreducible representations of
$\Gpo$
is done following the steps below:

1) One identifies the universal covering group
$\Gpocov$
of
$\Gpo$.
Let us denote by
$
\pi: \Gpocov \rightarrow \Gpo
$
the canonical projection. If
$V$ is a projective (unitary irreducible) representation of
$\Gpo$
then
\be
\tilde{V} \equiv V \circ \pi
\ee
is a projective (unitary irreducible) representation of
$\Gpocov$
verifying:
\be
\tilde{V}(g_{0}) = \lambda \times id
\label{lift}
\ee
for any
$
g_{0} \in Ker(\pi).
$
Here
$\lambda$
is a complex number of modulus 1.

It will be clear immediately why it is more convenient to classify,
up to unitary equivalence, the projective unitary irreducible
representations of
$\Gpocov$
verifying the condition (\ref{lift}).

2) One determines the cohomology group
$
H^{2}(\Gpocov,\R)
$
using infinitesimal arguments. This is possible because for
connected and simply connected Lie groups one has
$
H^{2}(G,\R) \simeq H^{2}(Lie(G),\R).
$

3) One selects a multiplier from every cohomology class and
tries to classify the unitary irreducible
$m$-representations of
$\Gpocov$
and afterwards selects those verifying the condition
(\ref{lift}). This can be done as in \cite{G2} generalizing
theorem 7.16 of \cite {V}. Alternatively, one could try to
connect the
$m$-representations
to some projective systems of imprimitivity.
This approach will be used in the following.

We will start with the Galilei group in 1+3 dimensions and
afterwards we will examine in detail the more interesting case
of 1+2 dimensions.

3.2 Let
$\Gpo$
the proper orthochronous Galilei group in 1+3 dimensions.

1) The universal covering group is \cite{V} set-theoretically:
$
SU(2) \times \R^{3} \times \R \times \R^{3}
$
with the composition law:
\be
(A_{1},{\bf v}_{1},\eta_{1},{\bf a}_{1}) \cdot
(A_{2},{\bf v}_{2},\eta_{2},{\bf a}_{2}) =
(A_{1}A_{2},{\bf v}_{1}+\delta(A_{1}){\bf v}_{2},\eta_{1}+\eta_{2},
{\bf a}_{1}+\delta(A_{1}){\bf a}_{2}+\eta_{2}{\bf v}_{1}).
\ee

Here
$
\delta:SU(2) \rightarrow SO(3)
$
is the covering map described by:
\be
[\delta(A)\cdot x] = A [x] A^{*}
\ee
where
\be
[x] \equiv x_{1} \sigma_{1} + x_{2} \sigma_{2} + x_{3} \sigma_{3}
\ee
$
(\sigma_{1}, \sigma_{2}, \sigma_{3}
$
being the Pauli matrices).
The covering map
$\delta$
extends trivially to a map (denoted also by
$\delta$)
from
$\Gpocov$
onto
$\Gpo$
which is the looked for covering map.

2) In \cite{V} theorem 7.42 describes the multiplier group of
$\Gpocov$.
Namely, every such multiplier is cohomologous to one of the
following form:
\be
m_{\tau}(r_{1},r_{2}) \equiv exp\left\{ i{\tau\over 2}\left[
{\bf a}_{1}\cdot \delta(A_{1}){\bf v}_{2} -
{\bf v}_{1}\cdot  \delta(A_{1}){\bf a}_{2} +
\eta_{2} {\bf v}_{1}\cdot \delta(A_{1}){\bf v}_{2}\right]\right\}.
\ee
Here
$
r_{i} = (A_{i},{\bf v}_{i},\eta_{i},{\bf a}_{i}),~(i = 1,2)
$
and
$
\tau \in \R.
$
Moreover
$
m_{\tau_{1}} \sim m_{\tau_{2}}
$
{\it iff}
$
\tau_{1} = \tau_{2}.
$

3) Let
$V$
be a (unitary)
$
m_{\tau}$-representation of
$\Gpo$.
We define:
\be
U_{A,{\bf v}} \equiv V_{A,{\bf v},0,{\bf 0}}
\ee
\be
W_{\eta,{\bf a}} \equiv V_{{\bf I},{\bf 0},\eta,{\bf a}}.
\ee

Then we have
\begin{prop}
(i)
$U$
and
$W$
are (unitary) representations of
$
G = SU(2) \times_{\delta} \R^{3}
$
and respectively
$
\R \times \R^{3}
$
i.e.
\be
U_{A_{1},{\bf v}_{1}}~U_{A_{2},{\bf v}_{2}} =
U_{A_{1}A_{2},{\bf v}_{1}+\delta(A_{1}){\bf v}_{2}}
\ee
\be
W_{\eta_{1},{\bf a}_{1}}~W_{\eta_{2},{\bf a}_{2}} =
W_{\eta_{1}+\eta_{2},{\bf a}_{1}+{\bf a}_{2}}.
\ee

One also has:
\be
U_{A,{\bf v}}~W_{\eta,{\bf a}}~U_{A,{\bf v}}^{-1} =
exp\left\{i\tau\left[{\bf v}\cdot\delta(A){\bf a}+
\eta{{\bf v}^{2}\over 2}\right]\right\}
W_{\eta,\delta(A){\bf a}+\eta{\bf v}}
\label{connection}
\ee

(ii) Conversely, if
$U$
and
$W$
are as above, let us define
\be
V_{A,{\bf v},\eta,{\bf a}} \equiv
exp\left\{i{\tau\over 2} {\bf a} \cdot {\bf v}\right\}
W_{\eta,{\bf a}}~U_{A,{\bf v}}.
\label{V}
\ee

Then
$V$ is a (unitary)
$
m_{\tau}$-representation of
$\Gpocov$.

(iii) Moreover, the commuting ring of
$V$
is isomorphic to the commuting ring of the pair
$
(U,W).
$
\end{prop}

We push the analysis a step further by considering the projector
valued measure
$P$
associated to
$W$:
\be
W_{\eta,{\bf a}} = \int e^{i(p_{0}\eta+{\bf p}\cdot {\bf a})}
dP(p_{0},{\bf p}).
\label{int}
\ee

Then one easily proves that (\ref{connection}) is equivalent to:
\be
U_{A,{\bf v}}~P_{E}~U_{A,{\bf v}}^{-1} =
P_{(A,{\bf v})\cdot E}
\ee
where the action of
$G$
on
$
\R \times \R^{2}
$
is:
\be
(A,{\bf v})\cdot [p_{0},{\bf p}] =\left [p_{0}+{\tau\over 2} {\bf v}^{2}-
{\bf v}\cdot\delta(A){\bf p},\delta(A){\bf p}-\tau {\bf v}\right].
\label{act}
\ee

So
$
(U,P)
$
is a (true) system of imprimitivity. Moreover the commuting ring
of
$V$ is isomorphic to the commuting ring of
$
(U,P).
$
To classify the unitary irreducible
$
m_{\tau}$-representations
of
$\Gpocov$,
one classifies the irreducible systems of imprimitivity for
$G$
based on
$
\R \times \R^{3}
$
relative to the action (\ref{act}). According to the lemma
\ref{ergo} they correspond to the
$G$-orbits
relative to this action. Then one reconstructs
$W$
from (\ref{int}) and uses (\ref{V}) to obtain the
$
m_{\tau}$-representation
$V.$
The classification of the systems of imprimitivity of
$G$
relative to (\ref{act}) is easy because we do not have the
complications due to the existence of the multipliers discussed
in detail in the preceeding chapter. We will consider only the
physically interesting case
$
\tau \not= 0.
$

The orbits of (\ref{act}) are:
$$
Z_{\rho} = \{[p_{0},{\bf p}]\vert {\bf p}^{2}-2\tau p_{0}=\rho\}~~
(\rho \in \R).
$$

So we study the systems of imprimitivity for
$
\Gpocov$
on
$
Z_{\rho}.
$
The stability subgroup of
$
\left[-{\rho\over 2\tau},{\bf 0}\right]
$
is:
$$
G_{0} = \{(A,{\bf 0})\vert A \in SU(2)\} \simeq SU(2).
$$

The unitary irreducible representations of
$G_{0}$
are, up to unitary equivalence, the well known
$
D^{(s)}
$
(with
$
s \in \N/2
$).
A convenient associated cocycle is simply:
\be
\phi^{(s)}((A,{\bf v}),[p_{0},{\bf p}]) = D^{(s)}(A).
\ee

If we identify
$
Z_{\rho} \simeq \R^{3}
$
via
$
\left[{1\over 2\tau}({\bf p}^{2}-\rho),{\bf p}\right]
\leftrightarrow {\bf p}
$
and consider on
$\R^{3}$
the Lebesgue
$\Gpocov$-invariant
measure
$
d{\bf p},
$
then it follows that the corresponding system of imprimitivity
acts in
$
{\cal H} = L^{2}(\R^{3},{\bf C}^{2s+1},d{\bf p})
$
according to:
\begin{eqnarray}
\left(U_{A,{\bf v}}~f\right)({\bf p}) & = &
D^{(s)}(A)~f(\delta(A)^{-1}({\bf p}+\tau{\bf v}))
\\
P_{E} & = & \chi_{E}.
\end{eqnarray}

One reconstructs
$W$
immediately:
\be
\left(W_{\eta,{\bf a}}~f\right)({\bf p}) =
exp\left\{i\left({{\bf p}^{2}-\rho\over 2\tau}\eta+{\bf p}\cdot
{\bf a}\right)\right\}~f({\bf p})
\label{Ab}
\ee
and,
using (\ref{V}) gets the expression of the
$
m_{\tau}$-representation as:
\begin{eqnarray}
\nonumber
\left(V_{A,{\bf v},\eta,{\bf a}}~f\right)({\bf p}) &=&
exp\left\{i\left({\tau\over 2} {\bf v}\cdot {\bf a}+
{{\bf p}^{2}\over 2\tau}\eta-{\rho\over 2\tau}\eta+
{\bf p}\cdot {\bf a}\right)\right\}
\\
&~&D^{(s)}(A)~f(\delta(A)^{-1}({\bf p}+\tau{\bf v}))
\label{rep}
\end{eqnarray}

For obvious reasons the factor
$
\left\{-i{\rho\over 2\tau}\eta\right\}
$
can be discarded and we reobtain the results of ch. IX,8 of
\cite{V}.

Finally, one notices that the condition (\ref{lift}) is
fulfilled. So, the unitary irreducible
$
m_{\tau}$-representations
of
$\Gpocov$
are indexed, beside
$
\tau \in R^{*}
$
by
$
s \in \N/2$,
are given by (\ref{rep}) (with
$
\rho = 0
$)
and induce projective unitary irreducible representations of
$\Gpo$.

3.3 Let
$\Gpo$
be the proper orthochronous Galilei group in 1+2 dimensions.

1) The covering group
$\Gpocov$
is the set
$
\R \times \R^{2} \times \R \times \R^{2}
$
with the composition law:
\begin{eqnarray}\nonumber
&~&(x_{1},{\bf v}_{1},\eta_{1},{\bf a}_{1}) \cdot
(x_{2},{\bf v}_{2},\eta_{2},{\bf a}_{2}) =
\\
&~&(x_{1}+x_{2},{\bf v}_{1}+\delta(x_{1}){\bf v}_{2},\eta_{1}+\eta_{2},
{\bf a}_{1}+\delta(x_{1}){\bf a}_{2}+\eta_{2}{\bf v}_{1}),
\end{eqnarray}
where
$
\delta: \R \rightarrow SO(2)
$
is the covering map described by:
\begin{eqnarray}
\delta(x) \equiv \left( \matrix{ cos(x) & sin(x) \cr -sin(x) &
cos(x) } \right).
\end{eqnarray}

The map
$\delta$
extends obviously to the covering map (denoted also by
$\delta$)
from
$\Gpocov$
onto
$\Gpo$.

2) The multiplier group for
$\Gpocov$
is described in \cite{G2}. One first identifies
$
Lie(\Gpocov) \simeq \R \times \R^{2}\times \R \times \R^{2}
$
with the Lie bracket:
\be
[(\alpha_{1},{\bf u}_{1},t_{1},{\bf x}_{1}),
(\alpha_{2},{\bf u}_{2},t_{2},{\bf x}_{2})] =$$
$$(0,A(\alpha_{1} {\bf u}_{2}-\alpha_{2} {\bf u}_{1}),0,
A(\alpha_{1} {\bf x}_{2}-\alpha_{2} {\bf x}_{1})+
t_{2} {\bf u}_{1}-t_{1} {\bf u}_{2}).
\ee

Here
\begin{eqnarray}\nonumber
A = \left( \matrix{0 & 1 \cr -1 & 0} \right).
\end{eqnarray}

Next, one computes by simple algebraic manipulations
$
H^{2}(Lie(\Gpocov),\R).
$

One finds out that every Lie algebraic cocycle
$\xi$
is cohomologous to one of the form
$
\tau \xi_{0} + F \xi_{1} + S \xi_{2}
$
where
\be
\xi_{0}((\alpha_{1},{\bf u}_{1},t_{1},{\bf x}_{1}),
(\alpha_{2},{\bf u}_{2},t_{2},{\bf x}_{2})) = {\bf x}_{1}\cdot
{\bf u}_{2} - {\bf x}_{1}\cdot {\bf u}_{2}
\ee
\be
\xi_{0}((\alpha_{1},{\bf u}_{1},t_{1},{\bf x}_{1}),
(\alpha_{2},{\bf u}_{2},t_{2},{\bf x}_{2})) =
\frac{1}{2}  <{\bf u}_{1},{\bf u}_{2}>
\ee
\be
\xi_{0}((\alpha_{1},{\bf u}_{1},t_{1},{\bf x}_{1}),
(\alpha_{2},{\bf u}_{2},t_{2},{\bf x}_{2})) =
\alpha_{1} t_{2} - \alpha_{2} t_{1}.
\ee

Here
$
<\cdot,\cdot>
$
is the sesquilinear form on
$
\R^{2}
$
given by:
\be
<{\bf u},{\bf v}> \equiv {\bf u}^{t} A {\bf v}.
\ee

Moreover
$
\tau_{1} \xi_{0} + F_{1} \xi_{1} + S_{1} \xi_{2} \sim
\tau_{2} \xi_{0} + F_{2} \xi_{1} + S_{2} \xi_{2}
$
{\it iff}
$
\tau_{1} = \tau_{2},~ F_{1} = F_{2}~,S_{1} = S_{2}.
$

Next, one determines
$
H^{2}(\Gpocov,\R)
$
and finds out that every multiplier of
$\Gpocov$
is cohomologous to one of the following form:
\begin{eqnarray}\nonumber
m_{\tau,F,S}(r_{1},r_{2}) &\equiv &
exp\left\{i{\tau\over2}
[{\bf a}_{1}\cdot \delta(x_{1}) {\bf v}_{2}
- {\bf v}_{1}\cdot \delta(x_{1}) {\bf a}_{2} + \eta_{2}
{\bf v}_{1}\cdot \delta(x_{1}) {\bf v}_{2}]\right\}
\\
&\times &exp\left\{{iF\over 2} < {\bf v}_{1},\delta(x_{1}) {\bf v}_{2} > +
iS \eta_{1} x_{2}\right\}
\end{eqnarray}

Here
$
r_{i} = (x_{i},{\bf v}_{i},\eta_{i},{\bf a}_{i})~(i = 1,2).
$
Moreover
$
m_{\tau_{1},F_{1},S_{1}} \sim m_{\tau_{2},F_{2},S_{2}}
$
{\it iff}
$
\tau_{1} = \tau_{2},~F_{1} = F_{2},~S_{1} = S_{2}.
$

3) Let
$V$
be a (unitary)
$
m_{\tau,F,S}$-representation of
$\Gpocov$.
We define
\be
U_{x,{\bf v}} \equiv V_{x,{\bf v},0,{\bf 0}}
\ee
\be
W_{\eta,{\bf a}} \equiv V_{0,{\bf 0},\eta,{\bf a}}.
\ee

Then we have

\begin{prop}
(i)
$U$
is a (unitary)
$
m_{F} \equiv m_{0,F,0}$-representation of
$
G \equiv \R \times_{\delta} \R^{2}:
$
\be
U_{x_{1},{\bf v}_{1}}~U_{x_{2},{\bf v}_{2}} =
exp\left\{ {iF\over 2} <{\bf v}_{1},\delta(x_{1}){\bf v}_{2}>
\right\} U_{x_{1}+x_{2},{\bf v}_{1}+\delta(x_{1}){\bf v}_{2}}.
\ee

$W$
is a (unitary) representation of the Abelian group
$
\R \times \R^{2}
$
i.e.
\be
W_{\eta_{1},{\bf a}_{1}}~W_{\eta_{2},{\bf a}_{2}} =
W_{\eta_{1}+\eta_{2},{\bf a}_{1}+{\bf a}_{2}}.
\ee
One also has:
\be
U_{x,{\bf v}}~W_{\eta,{\bf a}}~U_{x,{\bf v}}^{-1} =
exp\left\{ i\tau\left[{\bf v}\cdot \delta(x){\bf a}+\eta{{\bf v}^{2}
\over 2}\right]+iS\eta x\right\}~W_{\eta,\delta(x){\bf a}
+\eta{\bf v}}
\label{conn}
\ee

(ii) Conversely, if
$U$
and
$W$
are as above, let us define:
\be
V_{x,{\bf v},\eta,{\bf a}} \equiv exp\left\{ i{\tau\over 2}
{\bf a}\cdot {\bf v}+iS\eta x\right\}~
W_{\eta,{\bf a}}~U_{x,{\bf v}}.
\label{V'}
\ee

Then
$V$
is a (unitary)
$
m_{\tau,F,S}$-representation of
$\Gpocov$.

(iii) Moreover, the commuting ring of
$V$
is isomorphic to the commuting ring of the couple
$
(U,W).
$
\end{prop}

Like in subsection 3.2, we consider the projector valued measure
$P$
associated to
$W$:
\be
W_{\eta,{\bf a}} = \int e^{i(\eta p_{0}+{\bf a}\cdot {\bf p})}
dP(p_{0},{\bf p}).
\label{int'}
\ee

Then one easily proves that (\ref{conn}) is equivalent to:
\be
U_{x,{\bf v}}~P_{E}~U_{x,{\bf v}}^{-1} =
P_{(x,{\bf v})\cdot E}
\ee
where the action of
$G$
on
$
\R \times \R^{2}
$
is
\be
(x,{\bf v})\cdot [p_{0},{\bf p}] =
\left[p_{0}-{\bf v}\cdot \delta(x){\bf p}+{\tau
{\bf v}^{2}\over 2}-Sx,\delta(x){\bf p}-\tau{\bf v}\right].
\label{act'}
\ee

So
$
(U,P)
$
is a
$m_{F}$-system of imprimitivity. Moreover, the commuting ring of
$V$
is isomorphic to the commuting ring of
$
(U,P).
$
The classification of the unitary irreducible
$
m_{\tau,F,S}$-representations of
$\Gpocov$
is reduced to the classification of the irreducible
$
m_{F}$-systems of imprimitivity for
$G$
based on
$
\R \times \R^{2}
$
and relative to the action (\ref{act'}). Then
$W$
is obtained from (\ref{int'}) and
$V$
from (\ref{V'}).

The orbits of (\ref{act'}) are easy to compute.
We distinguish four cases which must be studied separatedly:

a) $\tau \not= 0,~S = 0$

$$
Z^{1}_{\rho} \equiv \{ [p_{0},{\bf p}] \vert
{\bf p}^{2}-2\tau p_{0} = \rho \};~\rho \in \R,
$$

b) $\tau \not= 0,~S \not= 0$

$$
Z^{2} \equiv \{ [p_{0},{\bf p}] \vert
{\bf p} \in \R^{2}, p_{0} \in \R \}
$$

c) $\tau = 0,~S = 0$

$$
Z^{3}_{r} \equiv \{ [p_{0},{\bf p}] \vert p_{0} \in \R,
{\bf p}^{2} = r^{2} \};~r \in \R_{+}
$$
$$
Z^{4}_{p_{0}} \equiv \{ [p_{0},{\bf 0}] \};~p_{0} \in \R
$$

d) $\tau = 0,~S \not= 0$

$$
Z^{5}_{r} \equiv \{ [p_{0},{\bf p}] \vert
{\bf p}^{2} = r^{2},p_{0} \in \R \}; r \in \R_{+} \cup \{0\}.
$$

We analyse them one by one:

a)
We have:
$$
G_{0} = G_{\left[-{\rho\over 2\tau},{\bf 0}\right]} =
\{(x,{\bf 0}) \vert x \in \R \}\simeq \R
$$

Because
$
m\vert_{G_{0} \times G_{0}} = 0
$
we have to consider only the unitary irreducible representations
of
$\R$.
They are
$
D^{(s)}~(s \in \R)
$
acting in {\bf C} as follows:

\be
D^{(s)}(x,{\bf 0}) = e^{isx}.
\ee

As an associated cocycle we can take:
\be
\phi^{(s)}((x,{\bf v}),[p_{0},{\bf p}]) = e^{isx}.
\ee

The corresponding system of imprimitivity acts in
$
{\cal H} = L^{2}(\R^{2},d{\bf p})
$
if one notices that it is possible to identify
$
Z^{1}_{\rho} \simeq \R^{2}
$
with the Lebesgue measure
$
d{\bf p},
$
like in subsection 3.2 One has:
\be
\left(U_{x,{\bf v}}~f\right) ({\bf p}) = e^{isx}
f(\delta(x)^{-1}({\bf p}+\tau{\bf v}))
\ee
\be
P_{E} = \chi_{E}.
\ee

We reconstruct
$W$
as in subsection 3.2 (see (\ref{Ab})) and then the representation
$V$.
One gets in this way the expression:
\begin{eqnarray}\nonumber
\left(V_{(x,{\bf v},\eta,{\bf a}}~f\right) ({\bf p}) &=&
exp\left\{i\left({\tau\over 2} {\bf a}\cdot{\bf v}+
{{\bf p}^{2}-\rho\over 2\tau}+{\bf a}\cdot{\bf p}\right)\right\}
\\
&\times &exp\left\{i\left({F\over 2\tau} <{\bf v},{\bf p}>+sx\right)\right\}
f(\delta(x)^{-1}({\bf p}+\tau{\bf v})
\end{eqnarray}

Like in 3.2 one can discard the factor
$
exp\left(-i\eta{\rho\over 2\tau}\right).
$
The condition (\ref{lift}) is fulfilled so we have obtained
a family of
$
m_{\tau,F,0}$-representations depending beside
$
\tau, F
$
on
$s \in \R$.

b)
In this case we have:
$$
G_{0} = G_{[0,{\bf 0}]} = \{(0,{\bf 0}) \}
$$
and it is clear that one has to consider only the trivial
representation of
$
G_{0}.$
The corresponding
$
m_{F}$-system
of imprimitivity acts in
$
L^{2}(\R \times \R^{2},dp_{0}\otimes d{\bf p})
$
according to:
\be
\left(U_{x,{\bf v}}~f\right)(p_{0},{\bf p}) =
f(p_{0}+{\bf v}\cdot {\bf p}+{1\over 2}\tau {\bf v}^{2}+Sx,
\delta(x)^{-1}({\bf p}+\tau {\bf v}))
\ee
\be
P_{E} = \chi_{E}.
\ee

The reconstruction of
$W$
is immediate:
\be
\left(W_{\eta,{\bf a}}~f\right)(p_{0},{\bf p}) =
e^{i(\eta p_{0}+{\bf a}\cdot {\bf p})}~f(p_{0},{\bf p})
\ee
and using (\ref{V'}), the expression of a
$
m_{\tau,F,S}$-representation
of
$\Gpocov$
is obtained:
\begin{eqnarray}\nonumber
\left(V_{x,{\bf v},\eta,{\bf a}} f\right)(p_{0},{\bf p}) &=&
exp\left\{i\left[{\tau\over 2} {\bf a}\cdot {\bf v}+
{\bf a}\cdot {\bf p}+{F\over 2\tau}<{\bf v},{\bf p}>+\eta(p_{0}+Sx)
\right]\right\}
\\
&~& f(p_{0}+{\bf v}\cdot {\bf p}+{1\over 2}\tau {\bf v}^{2}+Sx,
\delta(x)^{-1}({\bf p}+\tau{\bf v})).
\end{eqnarray}

The condition (\ref{lift}) is fulfilled so
$V$
induces a projective representation of
$\Gpo$.

c,d)
These cases correspond to
$
\tau = 0.
$
The orbite
$
Z^{3}_{r}
$
and
$
Z^{5}_{r}~(r \not= 0)
$
can be analysed simultaneously. We have:
$$
G_{0} = G_{[0,r{\bf e}_{1}]} = \left\{\left(2\pi
n,-{2\pi nS\over r}{\bf e_{1}}+\alpha {\bf e_{2}}\right) \vert
n \in {\bf Z}, \alpha \in \R \right\}
$$

It is convenient to denote
\be
(n,\alpha) \equiv \left(2\pi n,-{2\pi nS\over r}{\bf e}_{1}+
\alpha {\bf e}_{2}\right)
\ee
and to observe that:
\be
(n_{1},\alpha_{1})\cdot(n_{2},\alpha_{2}) =
(n_{1}+n_{2},\alpha_{1}+\alpha_{2})
\ee
i.e.
$
G_{0} \simeq {\bf Z} \times \R.
$
Next, one computes
$
m_{F}\vert_{G_{0} \times G_{0}}:
$
\be
m_{F}((n_{1},\alpha_{1}),(n_{2},\alpha_{2})) =
exp\left\{ {ik\over 2}(\alpha_{1} n_{2}-\alpha_{2} n_{1})\right\}
\ee
where
$
k \equiv {2\pi FS\over r}
$.
So, we have two subcases:

(i) $F = 0,~S \not= 0$ or $F \not= 0,~S = 0$.

In both situations we have
$
k = 0
$
so
$
m_{F}\vert_{G_{0} \times G_{0}} = 0.
$
It follws that we have to take into account true representations
of the group
$
G_{0}.
$
They are of the form
$
D^{(s,t)}~(s \in \R~(mod~1),~t \in \R)
$
and act in {\bf C} according to:
\be
D^{(s,t)}(n,\alpha) = e^{2\pi isn}~e^{it\alpha}.
\ee

An associated
$
m_{F}$-cocycle
can be computed using the formula (\ref{co}). The final
expression is rather complicated. Define the auxilliary function:
\be
\Phi((x,{\bf v}),[p_{0},{\bf p}]) \equiv \left\{ {iF\over 2r^{2}}
<\delta(x){\bf p},{\bf v}>\left[{\bf v}\cdot \delta(x){\bf p}-
2p_{0}+S(c_{\bf p}+c_{\delta(x){\bf p}}-x)\right]\right\}
\label{phi}
\ee
where
$
\R^{2} \ni {\bf p} \mapsto c_{\bf p} \in \R
$
is a Borel cross section verifying:
\be\nonumber
\delta(c_{\bf p}){\bf e}_{1} = {{\bf p}\over r}.
\ee

Then a cocycle associated to
$
D^{(s,t)}
$
can be taken to be:
\be
\phi^{(s,t)}((x,{\bf v}),[p_{0},{\bf p}]) =
\Phi((x,{\bf v}),[p_{0},{\bf p}]) exp\left\{i\left[{k\over r}
<\delta(x){\bf p},{\bf v}>+2\pi isx \right]\right\}.
\ee

The
$
m_{F}$-system of imprimitivity acts in
$
L^{2}(\R \times C_{r},dp_{0} \otimes d\Omega)
$
where
$
C_{r} \equiv \{{\bf p} \in \R^{2}\vert {\bf p}^{2} = r^{2}\}
$
and
$
d\Omega
$
is the Lebesgue measure on
$C_{r}.
$

We have:
\be
\left(U_{x,{\bf v}}~f\right)(p_{0},{\bf p}) =
\phi^{(s,t)}((x,{\bf v}),(x,{\bf v})^{-1}\cdot [p_{0},{\bf p}])~
f(p_{0}+{\bf v}\cdot {\bf p}+Sx,\delta(x)^{-1}{\bf p})
\ee
\be
P_{E} = \chi_{E}.
\ee

The corresponding expression for
$W$
is
\be
\left(W_{\eta,{\bf a}}~f\right)(p_{0},{\bf p}) =
e^{i(\eta p_{0}+{\bf a}\cdot {\bf p})}~f(p_{0},{\bf p}).
\ee

(ii) $F \not= 0,~S \not= 0$

In this situation
$
k \not= 0
$
and we have to classify some projective representations of
$
G_{0},
$
namely those representations
$
D_{n,\alpha}
$
verifying:
\be
D_{n_{1},\alpha_{1}}~D_{n_{2},\alpha_{2}} =
exp\{-ik(\alpha_{1} n_{2}-\alpha_{2} n_{1})\}~
D_{n_{1}+n_{2},\alpha_{1}+\alpha_{2}}
\label{D}
\ee

It is convenient to define:
\be
{\cal U}_{\alpha} \equiv D_{0,\alpha},~~~{\cal W}_{n} \equiv D_{n,0}.
\ee

Then we have:
\begin{prop}
(i)
$\U$
and
$\W$
are unitary representations of the Abelian groups
$\R$
and {\bf Z} respectively, i.e.:
\be
{\cal U}_{\alpha_{1}}~{\cal U}_{\alpha_{2}} =
{\cal U}_{\alpha_{1}+\alpha_{2}}
\ee
\be
{\cal W}_{n_{1}}~{\cal W}_{n_{2}} = {\cal W}_{n_{1}+n_{2}}.
\ee

We also have:
\be
{\cal U}_{\alpha}~{\cal W}_{n}~{\cal U}_{\alpha}^{-1} =
e^{-ik\alpha n}~{\cal W}_{n}.
\label{UW}
\ee

(ii) Conversely, if
$\U$
and
$\W$
are as above and we define
\be
D_{n,\alpha} \equiv e^{{i\over 2} kn\alpha}~
{\cal W}_{n}~{\cal U}_{\alpha}
\ee
then
$D$
is a projective unitary representation of
$
{\bf Z} \times \R
$;
more precisely it verifies (\ref{D}).

(iii) Moreover, the commuting ring of
$D$
is isomorphic to the commuting ring of the couple
$
(\U,\W).
$
\end{prop}

It is almost evident that the couple
$
(\U,\W)
$
can be connected to a system of imprimitivity applying SNAG
theorem. Indeed, the unitary irreducible representations of
{\bf Z} are one-dimensional and have the well-known expression:
\be\nonumber
n \mapsto z^{n},~~(z \in {\bf T})
\ee

So
$
\hat{\bf Z} \simeq {\bf T}
$
and SNAG theorem says that the representation
$\W$
has the generic form:
\be
{\cal W}_{n} = \int_{\bf T} z^{n}~dP(z)
\label{W}
\ee
where
$P$
is some projector valued measure on {\bf T}.

Then (\ref{UW}) is shown to be equivalet to:
\be
{\cal U}_{\alpha}~P_{E}~{\cal U}_{\alpha}^{-1} = P_{\alpha\cdot E}
\ee
where the action of
$\R$
on {\bf T} is:
\be
\alpha\cdot z \equiv e^{ik\alpha} z.
\label{a}
\ee

So
$
(\U,P)
$
is a system of imprimitivity for the Abelian group
$\R$
based on {\bf T} and relative to the action (\ref{a}). Moreover
$
(\U,P)
$
is irreducible {\it iff} the representation
$D$
we have started with is irreducible. To classify these irreducible
systems of imprimitivity is easy because for
$
k \not= 0
$
they are transitive. One has
$$
G_{0} = G_{\{1\}} = \left\{ {2\pi n\over k} \vert n \in {\bf Z} \right\}.
$$

The unitary irreducible representations of
$
G_{0}
$
are one-dimensional; they are indexed by
$
\lambda \in \R~(mod~1)
$
and have the expression:
\be
\pi^{(\lambda)} \left( {2\pi n\over k}\right) =
e^{2\pi i\lambda n}.
\ee

An associated cocycle is:
\be
\phi^{(\lambda)}(\alpha,z) = e^{ik\lambda\alpha}
\ee
and the corresponding system of imprimitivity acts in
$
{\cal K} = L^{2}({\bf T}, d\omega)~(d\omega\equiv
$
the Lebesgue measure on {\bf T}) according to:
\be
\left({\cal U}_{\alpha}~f\right)(z) = e^{ik\lambda\alpha}~
f(e^{-ik\alpha}z)
\ee
\be
P_{E} = \chi_{E}.
\ee

The expression of
$\W$
is obtained from (\ref{W}) and it is:
\be
\left({\cal W}_{n}~f\right)(z) = z^{n}~f(z).
\ee

Summing up, it follows that every projective unitary irreducible
representation of
$
{\bf Z} \times \R
$
verifying (\ref{D}) is unitary equivalent to one of the form
$
D^{(\lambda)}
$
which acts in
$\K$
according to
\be
\left(D^{(\lambda)}_{n,\alpha}~f\right)(z) =
e^{ik(\lambda-n/2)\alpha}~z^{n}~f(e^{-ik\alpha}z).
\ee

Moreover
$
D^{(\lambda_{1})} \simeq D^{(\lambda_{2})}
$
{\it iff}
$
\lambda_{1} - \lambda_{2} \in {\bf Z}.
$

A corresponding
$
m_{F}$-cocycle is:
\be\nonumber
\left(\phi^{(\lambda)}((x,{\bf v}),[p_{0},{\bf p}])~f\right)(z)
= \Phi((x,{\bf v}),[p_{0},{\bf p}])
\\
e^{2\pi ik\lambda x} z^{x}~f\left( e^{-{ik\over r}
<\delta(x){\bf p},{\bf v}>}z\right)
\ee
where
$\Phi$
has been defined before (see (\ref{phi})).

So, the looked for
$
m_{F}$-system
of imprimitivity
$
(U,P)
$
acts in
$
L^{2}(\R \times C_{r}, dp_{0} \otimes d\Omega,{\cal K}) \simeq
L^{2}(\R \times C_{r} \times {\bf T},dp_{0} \otimes d\Omega
\otimes d\omega)
$
according to
\begin{eqnarray}\nonumber
\left(U_{x,{\bf v}}~f\right)(p_{0},{\bf p},z) &=&
\Phi((x,{\bf v}),(x,{\bf v})^{-1}\cdot [p_{0},{\bf p}])
e^{-2\pi ik\lambda x}
\\
&~&z^{-x}~f\left(p_{0}+{\bf v}\cdot {\bf p}+Sx,
e^{-{ik\over r}<\delta(x){\bf p},{\bf v}>}z\right)
\end{eqnarray}
\be
\left(P_{E}~f\right)(p_{0},{\bf p},z) = \chi_{E}(p_{0},{\bf p})
f(p_{0},{\bf p},z).
\ee

The expression of
$W$
is:
\be
\left(W_{\eta,{\bf a}}~f\right)(p_{0},{\bf p},z) =
e^{i(\eta p_{0}+{\bf a}\cdot {\bf p})}~f(p_{0},{\bf p},z).
\ee

We still have to analyse the case (d) with
$
r = 0.
$
In this case:
$$
G_{0} = G_{[0,{\bf 0}]} = \{(0,{\bf v})\vert {\bf v} \in \R^{2} \}
\simeq \R^{2}.
$$

We have to classify all the unitary irreducible
$
m_{F}$-representations
of the Abelian group
$
\R^{2}
$
i.e.
\be
D_{{\bf v}_{1}}~D_{{\bf v}_{2}} = e^{-{iF\over 2}
<{\bf v}_{1},{\bf v}_{2}>}~D_{{\bf v}_{1}+{\bf v}_{2}}.
\ee

It is clear that
$D$
will be unitary equivalent to a unitary irreducible
representation of the canonical commutation relations. So
$
D \simeq \pi^{CCR} \equiv
$
the Schroedinger representation (see e.g. \cite{T}). The
corresponding system of imprimitivity acts in the Hilbert space
$
L^{2}(\R,dp_{0},\K^{CCR})
$
where
$
\K^{CCR}
$
is the representation space of
$
\pi^{CCR}.
$

We have explicitely:
\be
\left(U_{x,{\bf v}}~f\right)(p_{0}) = \pi^{CCR}_{\bf v}~f(p_{0}+Sx)
\ee
\be
P_{E} = \chi_{E}
\ee
which gives:
\be
\left(W_{\eta,{\bf a}}~f\right)(p_{0}) = e^{i\eta p_{0}}~f(p_{0}).
\ee

Finally, the trivial orbit
$
Z^{4}_{p_{0}}
$
corresponds to
$
m_{F}$-representations of
$G$.

\begin{rem}
The case
$
\tau \not= 0
$
can be simplified somehow if we make the following obsevation.
Let
$V$
be a
$
m_{\tau,F,S}$-representation.
We define:
\be
V'_{x,{\bf v},\eta,{\bf a}} \equiv
V_{x,{\bf v},\eta,{\bf a}+{F\over 2\tau}A{\bf v}}.
\ee

Then
$V'$
is a
$
m_{\tau,0,S}$-representation.
So in the case
$
\tau \not= 0
$
one can analyse only the situation
$
F = 0.
$
After the list of the irreducible unitary representations
$V'$
is obtained, one gets the corresponding
$V$
by making
$
{\bf a} \mapsto {\bf a}+{F\over 2\tau}A{\bf v}.
$

This remark appears in an infinitesimal form in \cite{B}.
\end{rem}

\end{document}